\documentstyle[12pt,page1]{article}
\newcommand{\sect}[1]{\section{#1}\setcounter{equation}{0}}

\font\mbn=msbm10 scaled \magstep1
\font\mbs=msbm7 scaled \magstep1
\font\mbss=msbm5 scaled \magstep1
\newfam\mbff
\textfont\mbff=\mbn
\scriptfont\mbff=\mbs
\scriptscriptfont\mbff=\mbss
\def\mbf{\fam\mbff}
\def\Re{{\mbf R}}

\def\Co{{\mbf C}}

\def\Di{{\mbf D}}
\newtheorem{Th}{Theorem}[section]

\newtheorem{C}[Th]{Corollary}

\newtheorem{Proposition}[Th]{Proposition}
\newtheorem{R}[Th]{Remark}

\author{Alexander Brudnyi\thanks
{Research supported in part by NSERC.
}\\
Department of Mathematics and Statistics\\
University of Calgary, Calgary\\
Alberta, Canada T2N 1N4}
\title{A Jensen Inequality for a Family of Analytic Functions}
\date{}
\begin{document}
\maketitle
\begin{abstract}
{We improve an estimate obtained in [Br] for the average number of
limit cycles of a planar polynomial vector field situated in a neighbourhood 
of the origin provided that the field in a larger neighbourhood is close 
enough to a linear center. The result follows from a new distributional
inequality for the number of zeros of a family of univariate holomorphic
functions depending holomorphically on a parameter.}
\end{abstract}
\sect{\hspace*{-1em}. Introduction.}
\noindent {\bf 1.1.} Let $f:=\{f_{v}\ ; v\in B(s,R)\}$, $R>1$, be a family of 
holomorphic (in the open unit disk $\Di_{1}$) functions depending 
holomorphically on a parameter $v$ varying in the open Euclidean ball 
$B(s,R)\subset\Co^{s}$. Hereafter $B(s,R)$ and $E(s,R)$
denote the open Euclidean balls in $\Co^{s}$ and $\Re^{s}$ centered at 0 of 
radius $R$, and $\Di_{t}:=B(1,t)\subset\Co$. Further, set
\begin{equation}\label{2}
{\cal N}_{f,\rho}(v):=\#\{z\in\overline{\Di}_{\rho};\ \ f_{v}(z)=0\};
\end{equation}
in addition, ${\cal N}_{f,\rho}(v)=+\infty$, if $f_{v}=0$ identically on 
$\Di_{1}$.  In [Br] we proved\\
{\bf Theorem.} {\em Assume that the family $f$ satisfies
$$
\frac{\sup_{v\in B(s,R)}\sup_{z\in\Di_{1}}|f_{v}(z)|}{\sup_{v\in B(s,1)}
\sup_{z\in\Di_{\rho}}|f_{v}(z)|}\leq M<\infty\ .
$$
Then for every $T\geq 0$ , the inequality
$$
|\{v\in E(s,1) \ ;\ {\cal N}_{f,\rho}(v)\geq T\}|\leq 
c_{1}se^{-c_{2}T/log M}\cdot |E(s,1)|
$$
holds with constants $c_{1}, c_{2}$ depending only on $\rho,R$ .
Here $|E|$ denotes the Lebesgue measure of $E\subset\Re^{s}$.}
\\
We applied this result to estimate the average number of isolated closed
trajectories (limit cycles) of a planar polynomial vector field situated in 
a neighbourhood of the origin provided that the field in a larger 
neighbourhood is close enough to a linear center. Central to this subject is
the second part of Hilbert's sixteenth problem asking whether the number of 
limit cycles of a planar polynomial vector 
field is always bounded in terms of its degree.
This is one of the few Hilbert  problems which remain unsolved
(for the main results see e.g. [B], [E], [FY], [I] and references there in).
According to Smale (see [S]) the global estimate
should be polynomial in the degree $d$ of the components of the vector field.
The mean estimate obtained in [Br] is of the order $(\log d)^{2}$ which is
substantially better than Smale's conjecture.
\begin{Th}{\rm ([Br,Th.A])}\label{te2}
Consider the equation 
\begin{equation}\label{poi}
\begin{array}{c}
\dot x=-y+F(x,y), \qquad  \dot y=\ x+G(x,y),\qquad  \qquad  \\
F(x,y)=\sum_{1\leq i+k\leq d}a_{ki}x^{k}y^{i},\ \ \ 
G(x,y)=\sum_{1\leq i+k\leq d}b_{ki}x^{k}y^{i}\ , 
\end{array}
\end{equation}
\begin{equation}\label{n}
\sum_{i,k}|a_{ki}|^{2}+\sum_{i,k}|b_{ki}|^{2}\leq N^{2}\ .
\end{equation}
Let $s:=d(d+3)$, $d\geq 2$, be the dimension of the space of real
coefficients of equation (\ref{poi}). Let $v$ be a point in this space,
and $C(v)$ be the number of limit cycles of the corresponding equation
(\ref{poi}) in the disk $D_{1/2}:=\{(x,y)\in\Re^{2};\
|x|^{2}+|y|^{2}\leq 1/4\}$. Inequality (\ref{n}) implies that
$v\in E(s,N)$. Assume that
\begin{equation}\label{n1}
N\leq\frac{1}{40\pi\sqrt{d}}\ .
\end{equation}
Then there exist absolute positive constants $C_{1},C_{2}$ such that
for any $T\geq 0$
$$
|\{v\in E(s,N);\ C(v)\geq T\}|\leq C_{1}se^{-C_{2}T/\log d}
\cdot |E(s,N)| \ .
$$
\end{Th}
As a corollary we obtain
\begin{C}\label{average}
Under assumptions of Theorem \ref{te2} the average number of limit cycles
of vector fields (\ref{poi}) in the disk $D_{1/2}$ is bounded 
from above by $c(\log d)^{2}$ with an absolute constant $c>0${\rm .}
\end{C}
The main purpose of the present paper is to improve the above estimates.

Let $F_{l}(x,y):=a_{10}x+a_{01}y$ and $G_{l}(x,y):=b_{10}x+b_{01}y$ be
the linear parts of $F$ and $G$. We
identify $\Re^{4}$ with the subspace of $\Re^{s}$ consisting of coefficients 
of the pair $F_{l},G_{l}$.
Let $\pi:\Re^{s}\longrightarrow\Re^{4}$ be the natural projection in the
space of coefficients and $V\subset\Re^{s}$ be a convex set such that
\begin{equation}\label{mes}
\frac{|E(4,N)|}{|\pi(V)|}\leq\delta<\infty\ .
\end{equation}
In what follows $|E|$ for $E\subset V$ denotes the Lebesgue measure of 
$E$ in $A$, where
$A\subset\Re^{s}$ is the real affine hull of $V$.

Under the assumptions of Theorem \ref{te2} and (\ref{mes}) we prove
\begin{Th}\label{te5}
There exist absolute positive constants $C_{1},C_{2},C_{3}$ such that
for any $T\geq 0$
$$
|\{v\in V;\ C(v)\geq T\}|\leq C_{1}\delta d^{C_{2}}e^{-C_{3}T}\cdot |V|
$$
\end{Th}
\begin{C}\label{average1}
There is an absolute positive constant $c>0$ such that
$$
\frac{1}{|V|}\int_{V}C(v)dv\leq c(\log\delta +\log d)\ .
$$
\end{C}
In particular, for $\delta=1$ we obtain that $(\log d)^{2}$ in Corollary 
\ref{average} can be replaced by $\log d$. In view of Smale's conjecture,
the above inequalities are especially interesting 
when $\delta$ is a polynomial function of $d$.
\\
{\bf 1.2.} We deduce Theorem \ref{te5} from a new
distributional inequality for families of univariate
holomorphic functions depending holomorphically on a parameter.
To formulate this result, we first recall some definitions from [Br1].

Let ${\cal O}_{R}$ denote the set of holomorphic functions defined on 
$B(N,R)\subset\Co^{N}$. In [Br1, Th.1.1] we proved the following statement.

{\em Let $f\in {\cal O}_{R}$, $R>1$, and $I$ be a real interval situated in
$B(N,1)$. (Hereafter we identify $\Co^{N}$ with $\Re^{2N}$.) There
is a constant $d=d(f,R)>0$ such that for any $I$ and any measurable subset
$\omega\subset I$ }
\begin{equation}\label{eq1}
\sup_{I}|f|\leq\left(\frac{4|I|}{|\omega|}\right)^{d}\sup_{\omega}|f|\ .
\end{equation}
The optimal constant in (\ref{eq1}) is called the {\em Chebyshev degree}
of $f\in {\cal O}_{R}$ in $B(N,1)$ and is denoted by $d_{f}(R)$. For 
instance, according to the classical Remez inequality $d_{f}(R)$ does not 
exceed the (total) degree of $f$, provided that $f$ is a polynomial. 
Example 1.14 in [Br1] shows that even in this case it can be essentially 
smaller than the degree. Note also that $d_{f}(R)$, can be estimated in terms 
of the valency of $f$ (see [Br1,Prop.1.7]).

We are ready to formulate our second result. 
\begin{Th}\label{te4}
Let $f:=\{f_{v}\in {\cal O}_{R}, R>1\}$ be a family of 
holomorphic in $\Di_{1}$ functions depending holomorphically on $v$. Define
$f_{0}(v):=f_{v}(0)\in {\cal O}_{R}$.
Let $V\subset B(s,1)\subset\Co^{s}(\cong\Re^{2s})$ be a convex set of 
real dimension $k$. Assume that
$$
\sup_{v\in V}\sup_{z\in\Di_{1}}|f_{v}(z)|\leq M_{1}<\infty\ \ \ 
{\rm and}\ \ \ \sup_{v\in V}|f_{0}(v)|\geq M_{2}>0 .
$$
Then for $\rho<1$ and every $T\geq 0$ , 
$$
|\{v\in V \ ;\ {\cal N}_{f,\rho}(v)\geq T\}|\leq 
4k\left(\frac{M_{1}}{M_{2}}\right)^{1/d_{f_{0}}(R)}\cdot 
e^{-[\log(1/\rho)/d_{f_{0}}(R)]T}\cdot |V|\ .
$$
\end{Th}
As a corollary we obtain
\begin{C}\label{c2}
Under the assumptions of Theorem \ref{te4},
$$
\frac{1}{|V|}\int_{V}{\cal N}_{f,\rho}(v)dv\leq
\frac{d_{f_{0}}(R)\log(4ek)+\log(M_{1}/M_{2})}{\log(1/\rho)}\ .
$$
\end{C}
\sect{\hspace*{-1em}. Proof of Theorem \ref{te4} and Corollary
\ref{c2}.}
\noindent {\bf Proof of Theorem \ref{te4}.}
We first recall the classical Jensen inequality (see e.g. [R,p.299]).

Let $f$ be a holomorphic function in $\Di_{1}$ continuous
on $S^{1}$ such that $f(0)\neq 0$. Let $n_{f}(r)$ be the
number of zeros of $f$ in $\overline{\Di}_{r}$, $0<r<1$, counted with their 
multiplicities, and
$M_{f}:=\sup_{S^{1}}|f|$. Then
\begin{equation}\label{iensen}
n_{f}(r)\leq\frac{\log(M_{f}/|f(0)|)}{\log(1/r)}\ .
\end{equation}

Let us denote 
$\omega(T):=\{v\in V \ ;\ {\cal N}_{f,\rho}(v)\geq T\}$. 
Then exactly as in [Br, Prop.2.1] one can 
check that
the function ${\cal N}_{f,\rho}$ is upper semicontinuous on $V\setminus S$, 
where $S$ is a certain closed subset of $V$ and $|S|=0$. In particular, 
$\omega(T)$ is measurable.
Without loss of generality we will assume that there is an $x\in\omega(T)$
such that $\sup_{\omega(T)}|f_{0}|=|f_{0}(x)|$. Then from
(\ref{iensen}) and [Br1,Th.1.9] we have
$$
T\leq {\cal N}_{f,\rho}(x)\leq\frac{\log(M_{1}/|f_{0}(x)|)}{\log(1/\rho)}
\ \ \ {\rm and} \ \ \
M_{2}\leq\sup_{V}|f_{0}|\leq
\left(\frac{4k |V|}{|\omega(T)|}\right)^{d_{f_{0}}(R)}
\cdot\sup_{\omega(T)}|f_{0}|\ .
$$
From these inequalities we obtain
$$
T\leq\frac{1}{\log(1/\rho)}\cdot\log\left(\frac{M_{1}}{M_{2}}\cdot
\left(\frac{4k|V|}{|\omega(T)|}\right)^{d_{f_{0}}(R)}\right)\ .
$$
The latter is equivalent to the required inequality.\ \ \ \ \ $\Box$\\
{\bf Proof of Corollary \ref{c2}.} Let 
$$
K:=\frac{\log[(M_{1}/M_{2})\cdot (4k)^{d_{f_{0}}(R)}]}{\log(1/\rho)}\ .
$$
Then a well-known formula and the inequality of 
Theorem \ref{te4} imply
$$
\begin{array}{c}
\displaystyle
\frac{1}{|V|}\int_{V}{\cal N}_{f,\rho}(v)dv=\int_{0}^{\infty}\omega(T)dT\leq\\
\\
\displaystyle
\int_{0}^{\infty}
\min\left\{1,\ 4k\left(\frac{M_{1}}{M_{2}}\right)^{1/d_{f_{0}}(R)}\cdot 
e^{-[\log(1/\rho)/d_{f_{0}}(R)]T}\right\}dT=\\
\\
\displaystyle
\int_{0}^{K}dT+
4k\left(\frac{M_{1}}{M_{2}}\right)^{1/d_{f_{0}}(R)}\int_{K}^{\infty}
e^{-[\log(1/\rho)/d_{f_{0}}(R)]T}dT=K+d_{f_{0}}(R)/\log(1/\rho)=\\
\\
\displaystyle
\frac{d_{f_{0}}(R)\log(4ek)+\log(M_{1}/M_{2})}{\log(1/\rho)}\ . \ \ \ \ \
\Box
\end{array}
$$
\sect{\hspace*{-1em}. Proof of Theorem \ref{te5} and Corollary 
\ref{average1}.}
\noindent 
For completness of the proof we repeat some arguments presented in [Br].

Passing to polar coordinates in equation (\ref{poi}) we get
\begin{equation}\label{zam}
\frac{dr}{d\phi}=\frac{P}{1+Q}r
\end{equation}
where $P(r,\phi):=\frac{xF(x,y)+yG(x,y)}{r^{2}},\
Q(r,\phi):=\frac{xG(x,y)-yF(x,y)}{r^{2}}$, 
$x=r\cos\phi,\ y=r\sin\phi$.
Let us complexify $r$: $r\in\Di_{1}=\{z\in\Co\ ;\ |z|<1\}$. Consider
equation (\ref{poi}) with complex coefficients that satisfy (\ref{n}).
In the domain $U=\Di_{1}\times [0,2\pi]$ we have
$$
\sup_{U}\left|\frac{F(x,y)}{r}\right|\leq\sqrt{\sum_{1\leq i+k\leq d}
|a_{ki}|^{2}}\sqrt{\sum_{1\leq i+k\leq d}(\cos\phi)^{2k}(\sin\phi)^{2i}}
\leq\sqrt{\sum_{1\leq i+k\leq d}|a_{ki}|^{2}}\cdot\sqrt{d}\ .
$$
Similarly,
$$
\sup_{U}\left|\frac{G(x,y)}{r}\right|\leq
\sqrt{\sum_{1\leq i+k\leq d}|b_{ki}|^{2}}\cdot\sqrt{d}\ .
$$
Hence, by (\ref{n}),
$$
\sup_{U}\left|\frac{P}{1+Q}\right|\leq\frac{N\sqrt{d}}{1-N\sqrt{d}}=:
\delta_{N}\ .
$$
For the $N$ of (\ref{n1}),

\begin{equation}\label{delt}
\delta_{2N}<3N\sqrt{d}\leq\frac{3}{40\pi}\ .
\end{equation}

\begin{Proposition}\label{approxim}
Consider the equation
\begin{equation}\label{Le}
\frac{dz}{d\phi}=H(z,\phi)\cdot z, \qquad (z,\phi)\in U,\qquad
\sup_{U}|H|\leq\delta_{2N},
\end{equation}
where $\delta_{2N}$ satisfies (\ref{delt}). Then any solution $z(\phi)$ of
(\ref{Le}) with initial condition $z(0)\in\Di_{3/4}$ may be extended to
$[0,2\pi]$, and $|z(\phi)-z(0)|\leq 8\pi N\sqrt{d}|z(0)|$ .
\end{Proposition}

{\bf Proof.} We have $\frac{d(\log z)}{d\phi}=H$. Hence, while
$|z(\phi)|\leq 1$, $\phi\in [0,2\pi]$, we have
$$
|z(\phi)|\leq e^{2\pi\delta_{2N}}|z(0)|\ .
$$
This follows from the Lagrange inequality. But $e^{2\pi\delta_{2N}}<
e^{\frac{3}{20}}<4/3$; hence, $|z(\phi)|\leq\frac{4}{3}|z(0)|< 1$.
In particular, $(z(\phi),\phi)\in U$ and
$$
z(\phi)=z(0)e^{\int_{0}^{\phi}H(z(t),t)dt}\ .
$$
Thus we have
$$
|z(\phi)-z(0)|\leq |z(0)|\cdot |e^{\int_{0}^{\phi}H(z(t),t)dt}-1|
\leq |z(0)| e^{2\pi\delta_{2N}}2\pi\delta_{2N}<8\pi N\sqrt{d}|z(0)|. 
$$

Now let $P_{v}$ be the Poincar\'{e} map corresponding
to equation
(\ref{zam}) obtained from a vector field $v$ of the type (\ref{poi}) with
complex coefficients. Consider the function
$g_{v}(z)=\frac{P_{v}(z)}{z}-1$ .
According to the method of successive approximations to solutions of
(\ref{zam}), $g_{v}(z)$ is holomorphic on 
$B(s,1/(20\pi\sqrt{d}))\times\Di_{3/4}$. The limit cycles of (\ref{poi}) 
located in
$D_{1/2}$ correspond to certain zeros of $g_{v}(z)$ in $\Di_{1/2}$. For
$|z|< 3/4$ and $v\in B(s,2N)$ we have $|g_{v}(z)|\leq 8\pi N\sqrt{d}$, by
Proposition \ref{approxim}. 

Let us consider linearization of (\ref{poi}) 
\begin{equation}\label{li}
\dot x=-y+F_{l}(x,y), \qquad  \dot y=\ x+G_{l}(x,y)\ .
\end{equation}
Let $w:=(a_{10},a_{01},b_{10},b_{01})\in\Co^{4}$ be a point 
in the space of coefficients of (\ref{li}). The
Poincar\'{e} map for (\ref{li}) can be calculated explicitly:
\begin{equation}\label{linear}
\begin{array}{c}
\displaystyle
P_{w}^{l}(z):=e^{f(w)}\cdot z;\\
\\
\displaystyle
f(w):=\int_{0}^{2\pi}
\frac{a_{10}\cos^{2}\phi+b_{01}\sin^{2}\phi+(a_{01}+b_{10})\sin\phi
\cdot\cos\phi}{1+b_{10}\cos^{2}\phi-a_{01}\sin^{2}\phi+(b_{01}-a_{10})\sin\phi
\cdot\cos\phi}d\phi\ .
\end{array}
\end{equation}
Also it is easy to see that 
$$
g_{0}(v):=g_{v}(0)=e^{f(\pi(v))}-1\ .
$$
Let $v_{0}$ be the vector field defined by
$$
\dot x=-y+\frac{N}{\sqrt{2}}x,\ \ \ \dot y=x+\frac{N}{\sqrt{2}}y\ .
$$
By definition, $v_{0}\in E(s,N)$ and 
$$
g_{0}(v_{0})=e^{\sqrt{2}\pi N}-1>\sqrt{2}\pi N\ .
$$
Further, $h(v):=g_{0}(v/(40\pi\sqrt{d}))\in {\cal O}_{2}$.
Let $D:=d_{h}(2)$ be the Chebyshev degree of $h$ in the ball $B(s,1)$.
Since $h$ is pullback by $\pi$ of a holomorphic function defined
on $B(4,2)\subset\Co^{4}$, $D$ coincides with
the Chebyshev degree of $h|_{B(4,2)}$ in $B(4,1)$ and therefore it does not 
depend on $s$. 

Now assume that $V\subset E(4,N)$ satisfies (\ref{mes}).
Then [Br1,Th.1.9] applied to $g_{0}|_{E(4,N)}$ implies that
$$
\sup_{\pi(V)}|g_{0}|_{E(4,N)}|
\geq (1/16\delta)^{D}\sup_{E(4,N)}|g_{0}|_{E(4,N)}|
\geq (1/16\delta)^{D}|g_{0}(v_{0})|>\sqrt{2}\pi N(1/16\delta)^{D}\ .
$$
But $\sup_{\pi(V)}|g_{0}|_{E(4,N)}|=\sup_{V}|g_{0}|$.
Thus applying Theorem \ref{te4} to $g(v,z):=g_{v}(3z/4)$ 
with $M_{1}=8\pi N\sqrt{d}$ and $M_{2}=\sqrt{2}\pi N(1/16\delta)^{D}$,
and using the fact that $k\leq s<3d^{2}$ (for $d\geq 2$), we have
$$
|\{v\in V;\ C(v)\geq T\}|\leq|\{v\in V \ ;\ {\cal N}_{g,2/3}(v)\geq T\}|\leq 
C_{1}\delta d^{C_{2}}\cdot 
e^{-C_{3}T}\cdot |V|\ 
$$
with $C_{1}=6\cdot 32^{1+1/2D}$, $C_{2}=2+1/2D$ and $C_{3}=\log(3/2)/D$.
\ \ \ \ \ $\Box$\\
\begin{R}\label{corr}
{\rm We used in the proof the fact that $\{v\in V;\ C(v)\geq T\}$ is a 
measurable
set. Indeed, $C(v)$ coincides with the number of nonnegative zeros of 
$g_{v}(z)$. Then
a slight modification of the proof of Proposition 2.1 in [Br] shows
that outside of a set of measure 0 in $V$, the function $C(v)$ is the
pointwise limit of a nonincreasing sequence of measurable
functions. In particular, $C(v)$ is measurable.}
\end{R}
{\bf Proof of Corollary \ref{average1}.} The result follows directly from
Corollary \ref{c2}. \ \ \ \ \ $\Box$
\begin{R}\label{concl}
{\rm In general the inequality of [Br,Th.B] (see the beginning of this paper)
is more powerful than that of Theorem \ref{te4}.
\\
The methods presented in this paper and in [Br] can be applied for 
more general equations of type (\ref{poi}), where $F$, $G$ are analytic 
functions of $x$ and $y$ depending analytically on a multidimensional 
parameter.}
\end{R}

\end{document}